\documentclass[11pt]{amsart}
\usepackage{graphicx}
\usepackage{amsfonts}
\usepackage{amscd}
\usepackage{amssymb}
\usepackage{alltt}

\newcommand{\ring}[1]{\mathbb{#1}}

\def\HH{{\mathcal H}}

\def\fourth{{\root 4 \of 12}}
\def\fpiZ{4\pi\ring{Z}}
\def\leftopen{(}
\def\rightopen{)}
\def\leftclosed{[}
\def\rightclosed{]}

\newcommand{\be}{\begin{equation}}
\newcommand{\ee}{\end{equation}}

\def\1{{\mu\mkern-6mu\mu}}

\def\op#1{{\operatorname{#1}}}

\title{The Honeycomb Problem on the Sphere}
\author{Thomas C. Hales}

\date{November 15, 2002}

\begin{document}

\theoremstyle{plain}
\newtheorem{thm}{Theorem}
\newtheorem{lemma}[thm]{Lemma}
\newtheorem{cor}[thm]{Corollary}
\newtheorem{prop}[thm]{Proposition}

\theoremstyle{definition}
\newtheorem{rem}[thm]{Remark}
\newtheorem{defn}[thm]{Definition}
\newtheorem{ex}[thm]{Example}
\newtheorem{estimate}[thm]{Estimate}

\begin{abstract}  The honeycomb problem on the sphere asks for
the partition of the sphere into $N$ equal areas that minimizes the
perimeter. This article solves the problem when $N=12$.  The unique
minimizer is a tiling of $12$ regular pentagons in the dodecahedral
arrangement.
\end{abstract}
\maketitle

\section{Introduction}

    The classical honeycomb problem asks for the perimeter-minimizing partition
of the plane into regions of equal area.  The optimal solution is the
regular hexagonal honeycomb tiling.  This article adapts the planar
proof to partitions on the unit sphere.  The honeycomb problem on the
sphere asks for the perimeter-minimizing partition of the sphere into
$N$ equal areas.  This article solves the problem when $N=12$.  The
unique minimizer is a tiling of $12$ regular pentagons in the
dodecahedral arrangements.  This problem was solved by L. Fejes T\'oth,
under an assumption of convexity in \cite{FT43}.

\section{The hexagonal isoperimetric inequality}
\label{sec:hex}

    The honeycomb theorem in the plane follows by showing that the regular
hexagon is the unique minimizer of a particular functional on the set of
closed curves with finitely many marked points.  This section gives a
review of this result.  It is the result that will be generalized
in this article.

\begin{figure}[htb]
  \centering
  \includegraphics{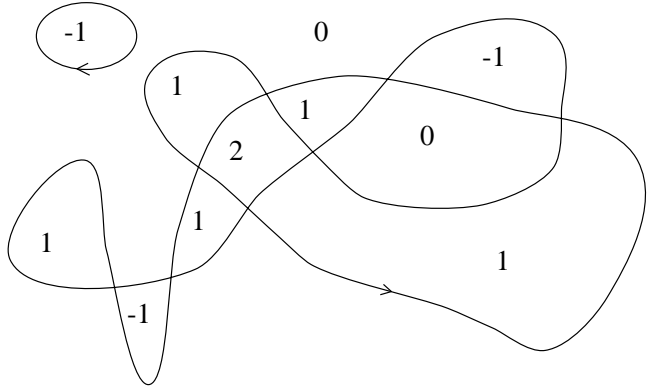}
  \caption{Signed areas in the plane.}
  \label{diag:signed}
\end{figure}

Let $\Gamma$ be a closed piecewise simple rectifiable curve in the
plane. We use the parameterization of the curve to give it a direction,
and use the direction to assign a signed area to the bounded components
of the plane determined by the curve. For example, if $\Gamma$ is a
piecewise smooth curve, the signed area is given by Green's formula
    $$\int_{\Gamma} x dy.$$
Generally, we view $\Gamma$ as an integral current \cite[p.44]{M95}. We
let $P$ be an integral current with boundary $\Gamma$. Expressed
differently, we give a signed area by assigning an multiplicity
    $m(U)\in\ring{Z}$
to each bounded component $U$ of $\ring{R}^2\setminus\Gamma$. (An
illustration appears in Figure~\ref{diag:signed}.) The area is $\sum
m(U)\op{area}(U)$. $P$ is represented by the formal sum
    $P = \sum m(U) U$.

Let $v_1,\ldots,v_t$, $t\ge 2$, be a finite list of points on $\Gamma$.
We do not assume that the points are distinct. Index the points
$v_{1},v_{2},\ldots,v_{t}$, in the order provided by the
parameterization of $\Gamma_i$.  Join $v_{i}$ to $v_{i+1}$ by a directed
line segment $f_{i}$ (take $v_{t+1}=v_{1}$). The chords $f_i$ form a
generalized polygon, and from the direction assigned to the edges, it
has a signed area $A_P\in\ring{R}$.

Let $e_i$ be the segment of $\Gamma$ between $v_{i}$ to $v_{i+1}$.  Let
$f^{\op{op}}$ be the chord $f$ with the orientation reversed. Let
$x(e_{i})\in\ring{R}$ be the signed area of the integral current bounded
by $(e_{i},f_{i}^{\op{op}})$. Let $E(P)=\{e_{i}\}$ denote the set of
edges of $P$. Accounting for multiplicities and orientations, we have
$$\op{area}(P) = A_P + \sum_{e\in E(P)} x(e).$$
Let $\alpha(P) =\min(1,\op{area}(P))$.

Define a truncation function $\tau:\ring{R}\to \ring{R}$ by
$$\tau(x) = \begin{cases}
        1/2, & x\ge 1/2,\\
        x,   & |x|\le 1/2,\\
        -1/2, &x\le -1/2.\end{cases}
$$
Set $T(P)= \sum_{E(P)}\tau(x(e))$. Recall that the perimeter of a
regular hexagon of unit area is $2\fourth$. Let $L(P)$ be the length of
$\Gamma$. Let $N(P)$ be the number of points $v_i$ on $\Gamma$, counted
with multiplicities. Let $a(N) = \min(2\pi\sqrt3/(3N^2),1)$.  The
following is proved in \cite{H}.

\begin{thm} (hexagonal isoperimetric inequality)
\label{thm:hex} %
Define $P$, $L(P)$, $N(P)$, and $a(N)$ as above. Assume that the signed
area of $P$ is at least $a(N(P))$. Then
    $$L(P)\ge - T(P)\fourth - (N(P)-6) 0.0505 + 2\alpha(P)\fourth.$$
Equality is attained if and only if $P$ is a regular hexagon
of area 1.
\end{thm}

\section{The pentagonal isoperimetric inequality}
\label{sec:pent}

We adapt the inequality of Theorem~\ref{thm:hex} to pentagons on the
unit sphere.

The first important difference as we move from the plane to the sphere
is that a simple closed curve bounds two simply connected regions. A
parameterized simple closed curved can be viewed as giving a positive
orientation to one of the regions and a negative orientation to the
other. It is not generally clear which area is to be preferred. This
ambiguity disappears if we take all areas on the sphere to be defined
modulo multiples of $4\pi$.  Thus, we define the area bounded by a
closed curve as taking values in $\ring{R}/(\fpiZ)$.  If a parameterized
curve bounds an area $A\in\ring{R}/(4\pi\ring{Z})$, then the oppositely
parameterized curve bounds an area $-A$.  A representative in $\ring{R}$
for a constant in $\ring{R}/(4\pi\ring{Z})$ will be called a {\it real
representative}. We define $\op{mrr}:
\ring{R}/(4\pi\ring{Z})\to\ring{Z}$ (the minimal real representative) to
be the function taking each $x$ to the real representative in
$\leftclosed -2\pi,2\pi\rightopen$.

We will also make use of the unsigned area of regions, which is the
usual Lebesgue measure of a region taking values in the set of
non-negative real numbers.  Write $\op{unsigned}(R)$ for the unsigned
area of a region $R$.

Let $\Gamma$ be a closed
piecewise simple rectifiable curve on the unit sphere.

As before, let $v_1,\ldots,v_t$, $t\ge 2$, be a finite list of points on
$\Gamma$. We do not assume that the points are distinct. Index the
points $v_{1},v_{2},\ldots,v_{t}$, in the order provided by the
parameterization of $\Gamma$.  Join $v_{i}$ to $v_{i+1}$ by an arc
$f_{i}$ of a great circle. (The choice of arc will be fixed in a
moment.) Take $v_{t+1}=v_{1}$. Let $P$ be an integral current with
boundary $\Gamma$. The chords $f_i$ form a generalized polygon, and from
the direction assigned to the edges, it has a signed area
$A_P\in\ring{R}/(\fpiZ)$.

\begin{figure}[htb]
  \centering
  \includegraphics{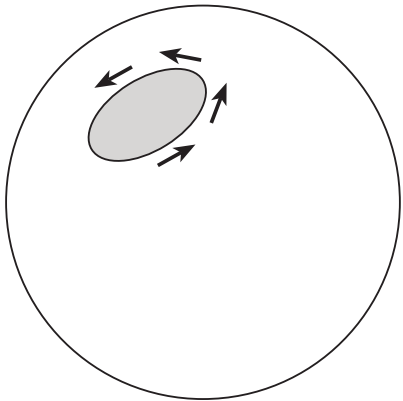}
  \caption{Signed areas on the unit sphere.  The two regions have
    signed areas $\epsilon$ and $-4\pi+\epsilon$}
  \label{diag:signed:sphere}
\end{figure}

Let $e_i$ be segment of $\Gamma$ between $v_{i}$
to $v_{i+1}$.  Let $f^{\op{op}}$ be the chord $f$ with the
orientation reversed.
Let $x(e_{i})\in\ring{R}/(\fpiZ)$ be the signed area of the integral current
bounded by $(e_{i},f_{i}^{\op{op}})$.

If $v_i$ and $v_{i+1}$ are antipodal, there are infinitely many geodesic
arcs $f_i$ joining $v_i$ and $v_{i+1}$.  In this case, choose $f_i$ in
such a way that $x(e_i)=0$.   If $v_i$ and $v_{i+1}$ are not antipodal,
then these two points break the great circle through them into two
circular arcs $f'$ and $f''$.  There is a difference of $2\pi$ between
the signed area of the integral current bounded by $(e_i,f'^{\op{op}})$
and that bounded by $(e_i,f''^{\op{op}})$.  We choose $f_i= f'$ or $f''$
in a way to make $x(e_{i})$ have its minimal real representative in
$\leftclosed -\pi,\pi \rightopen$.

Let $E(P)=\{e_{i}\}$ denote the set of edges of $P$. Accounting for
multiplicities and orientations, we have
    $$\op{area}(P) = A_P + \sum_{e\in E(P)} x(e)\in \ring{R}/(\fpiZ).$$
We emphasize that this identity does not hold in general between the
minimal real representatives of the terms.

Define a truncation function $\tau:\ring{R}/(\fpiZ)\to \ring{R}$ by
$$\tau(x) = \begin{cases}
        0.32, & \op{mrr}(x)\ge 0.32,\\
        \op{mrr}(x),   & |\op{mrr}(x)|\le 0.32,\\
        -0.32, &\op{mrr}(x)\le -0.32\end{cases}
$$
Set $\tau_0 = 0.32$.  Set $T(P)= \sum_{E(P)}\tau(x(e))$.

Let $\op{reg-perim}: \ring{R}/(\fpiZ)\times \ring{N}\to \ring{R}$ be the
function that gives the perimeter of a regular spherical $n$-gon of area
$x$.  It satisfies
    $$\op{reg-perim}(-x,n) = \op{reg-perim}(x,n).$$
The
perimeter of a regular spherical pentagon in the tiling by $12$
pentagons is
    $$\op{reg-perim}(\pi/3,5)\approx 3.64864.$$
Let $p_5$ equal this constant.

The function $\op{reg-perim}$ is given explicitly by the following
Mathematica module.
\bigskip   
\begin{alltt}
regPerim[area_, n_] := Module[\{cosgamma, alpha\},
    alpha = Pi - (2Pi - area)/n;
    cosgamma = (Cos[2Pi/n] + Cos[alpha/2]^2)/Sin[alpha/2]^2;
    n ArcCos[cosgamma]
    ]
\end{alltt}
\bigskip

\begin{figure}[htb]
  \centering
  \includegraphics{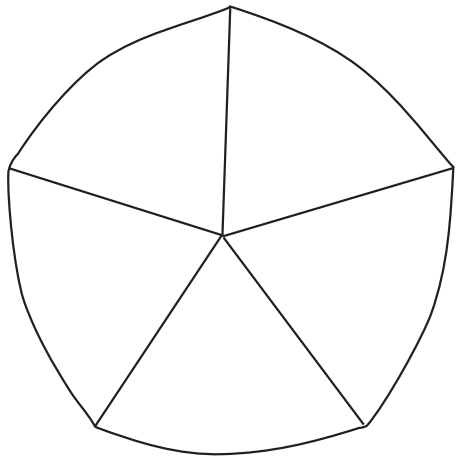}
  \caption{The triangulation used to compute $\op{reg-perim}$.}
  \label{diag:reg-perim}
\end{figure}

The derivation is a simple calculation based on spherical trigonometry
(the spherical law of cosines) and Figure~\ref{diag:reg-perim}. From
this explicit formula, we see that $\op{reg-perim}(\pi/3,n)$ extends to
an analytic function of $n>0$.  Let $p'_5$ be the partial derivative
$\partial_2\op{reg-perim}(\pi/3,5)$.

Let $\op{circ-perim}: \leftclosed 0,2\pi\rightopen\times
    \ring{R}/(\fpiZ)\to\ring{R}$ be defined as follows.
Let $\ell\in\leftclosed 0,2\pi\rightopen$ be the length of an arc of a
great circle on the unit sphere.  Draw a circle on the unit sphere
passing through the two endpoints of the arc with the property that one
of the two areas bounded by the arc and the circle has signed area $x$.
Let $\op{circ-perim}(\ell,x)$ be the length of the part of the circle
between the two endpoints that (together with the arc) bounds the region
of area $x$.  For example,
    $$\op{circ-perim}(0,x)$$
is the perimeter of a circle of area $x$.  Explicitly,
    \begin{equation}
    \op{circ-perim}(0,x) = \sqrt{4\pi x - x^2}.
    \label{eqn:circ-perim}
    \end{equation}
(Compare \cite[Eqn 10.1]{MAMM}.) To give another example,
    $$\op{circ-perim}(\pi,x) = \pi,$$
because in this case, the two endpoints of the arc are antipodal and the
circle must be a great circle through the antipodal points.

\begin{figure}[htb]
  \centering
  \includegraphics{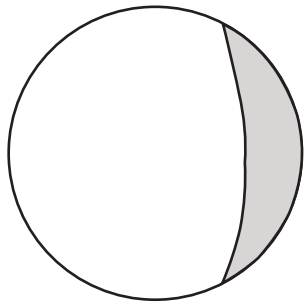}
  \caption{The arc of the circle bounding the
    shaded region has length $\op{circ-perim}$.}
  \label{diag:circ-perim}
\end{figure}

\begin{figure}[htb]
  \centering
  \includegraphics{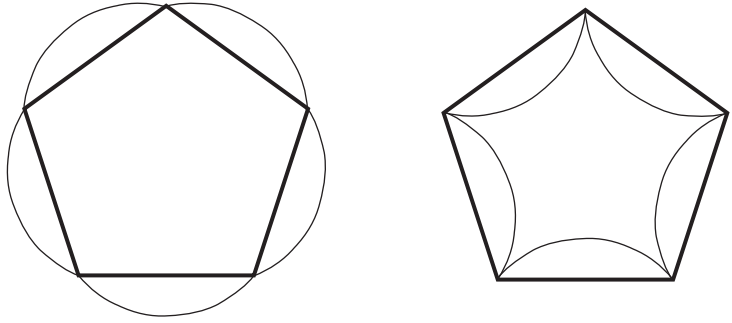}
  \caption{$\op{pent}(x)$ is the perimeter of bulging pentagon.
    The heavy set geodesic pentagon has area $\pi/3-x$ and each of
    the five two-sided regions has area $|x|/5$.}
  \label{diag:bulge}
\end{figure}

Let
    $$\op{pent}(x) = \op{circ-perim}(\op{reg-perim}
        (\pi/3 - x,5)/5,|x|/5) 5$$
It equals the perimeter of a regular pentagon of area $\pi/3$ in which
the geodesic edges have been replaced by arcs of a circle.  When $x>0$ the
edges bulge outward and when $x<0$, the edges bulge inward.
See Figure~\ref{diag:bulge}.  Because of the
absolute value in the formula, it is not obvious that this is differentiable
at $x=0$.  Nevertheless, it turns out to be, and we set
    $$B' = \op{pent}'(0).$$
A calculation shows in fact that
    $$B'=-\partial_1\op{reg-perim}_1(\pi/3,5)\approx -1.51.$$

Let $L(P)$ be the length of $\Gamma$. Let $N(P)$ be the number of points
$v_i$ on $\Gamma$, counted with multiplicities. Let
    $$a(N) =\min(3.75/n^2,0.1).$$
Let
    $$\rho(P) =\begin{cases}
                3\,|\op{mrr}(\op{mrr}(P))|/\pi&\text{if }
                    0\le|\op{mrr}(\op{area}(P))|\le \pi/3,\\
                1 & \text{otherwise.}
            \end{cases}
    $$

\begin{thm} (Isoperimetric inequality for spherical pentagons)
\label{thm:isoperimetric} Define $P$, $L(P)$, $N(P)$,  $a(N)$,
$\rho(N)$, $p_5$, $p_5'$, and $B'$ as above. Assume that
    $$\op{mrr}(\op{area}(P))\in\leftopen a(N(P)), 2\pi\rightopen.$$
Then
$$L(P)\ge T(P) B' - (5-N(P)) p_5' + \rho(P) p_5.$$
Equality is attained if and only if $P$ is a regular spherical pentagon
of unsigned area $\pi/3$, so that $L(P)=p_5$, $T(P)=0$, $N(P)=5$, and
$\rho(P)=1$.
\end{thm}

\section{The proof of the main result}

This section assumes Theorem~\ref{thm:isoperimetric} and derives the
main result from it.  The following section will give a proof of the
isoperimetric inequality.

\begin{thm} ($12$ honeycombs on a sphere)
Let $R_1,\ldots,R_{12}$ be a partition of the unit sphere into $12$
measurable sets of equal area $(\pi/3)$. Let $h = \HH^1(\cup
\{R_1,\ldots,R_{12}\})$ be the $1$-dimensional Hausdorff measure of the
current boundary.    Then
    $$ h \ge 6\,\op{reg-perim}(\pi/3,5).$$
If equality is attained, then up to a set of $1$-dimensional Hausdorff
measure zero, the current boundary is equal to the perimeter of the tiling
by $12$ congruent regular spherical pentagons.
\end{thm}

In preparation for the proof, we let $R_1,\ldots,R_{12}$ be a minimizer
of the $1$-dimensional Hausdorff boundary subject to the constraints
that each area is $\pi/3$.  A minimizer exists by \cite{M94}.  The
theorem is stated for planar soap bubbles; but the proof carries through
without difficulty to a sphere.

By the results of \cite{M94}, we may assume various regularity results.
Each $R_i$ is an open region bounded by finitely many analytic arcs. The
collection of arcs meet one another only at their endpoints. The number
of arcs meeting at each endpoint is two or three. Call this the {\it
degree\/} of the endpoint. Each $R_i$ consists of finitely many
connected components and each connected component is connected and
simply connected. Each connected component $P$ has a finite number of
marked points $N(P)$, determined by the endpoints of arcs of degree
three. Each connected component has at least one marked point.

Orient each region by the outward normal on the sphere, so that each
{\it signed area\/} $\op{area}(P)$ is positive and equal to the Lebesgue
measure of $R_i$.

Since the dodecahedral tiling has perimeter $6\,\op{reg-perim}(\pi/3,5)$,
the minimizer has a perimeter no greater than this constant.

\begin{lemma} In an optimal configuration,
the total unsigned area of  the connected components of unsigned area at
most $0.1$ is less than $4.0$.
\end{lemma}

\begin{proof}  Each connected component of unsigned area
    $x\in\leftopen0,0.1\rightclosed$
gives a perimeter at least $\op{circ-perim}(0,x)$. The function has
negative second derivative so that the perimeter is at least
    $$x \,\op{circ-perim}(0,0.1)/0.1.$$
Each segment of perimeter is shared by at most two such regions.  Thus,
if the total area is at least $4.0$, the total perimeter is at least
    $$4.0\,\op{circ-perim}(0,0.1)/0.2 > 6\, \op{reg-perim}(\pi/3,5),$$
contrary to the supposed optimality of the perimeter.
\end{proof}

Assume for a contradiction that we have found a strict contradiction to the
theorem, that is, a partition into equal areas with perimeter strictly
less than $6\,\op{reg-perim}(\pi/3,5)$.  Write this in the form
    \begin{equation}
    \sum_i (\op{perim}(R'_i) - \rho_i p_5) < 0,
    \label{eqn:CE}
    \end{equation}
where the sum extends over connected components $R'_i$, and the area
fractions $\rho_i$ are at most $1$ and sum to $12$.

We repeatedly modify the set of connected components by picking a
connected component of unsigned area at most $a(N)$, and deleting the
edge segments shared with another region. This decreases the number of
connected components by one.  We continue until all regions satisfy the
lower bound $a(N)$. The regions before this process had areas at most
$\pi/3$.  After the process, because of the lemma, each region $R_i''$
has area less than $\pi/3+4.0<2\pi$.

In the previous paragraph we always delete the edge segments with the
region that causes the perimeter to decrease the most.  Some edge
segment has length at least $1/N$ times the perimeter, where $N$ is the
number of junctures of degree $3$.  The perimeter, by the isoperimetric
inequality for spheres \cite{MAMM}, is at most $\op{circ-perim}(0,x)$
(that is, the perimeter of a circle of area $x$). Thus, the perimeter
decreases by at least $\op{circ-perim}(0,x)/N$, where $x =
\rho_i\pi/3\le a(N)$.

The sum $\sum_i \rho_i$ is unaffected by merging the two regions (say
with indices $1$ and $2$) if $\rho_1+\rho_2\le1$; that is, $\rho'_1 =
\rho_1+\rho_2$, $\rho_2'=0$.  But the sum can decrease if
$\rho_1+\rho_2>1$; that is, $\rho'_1 = 1$, $\rho_2'=0$, and
$\rho_1+\rho_2\ne \rho_1'+\rho_2'$. Overall, the left-hand side of
Inequality~\ref{eqn:CE} is at most what is was originally plus the term
    $$-2\op{circ-perim}(0,\rho_i\pi/3)/N + \rho_i p_5
    \le -2\op{circ-perim}(0,|x|)/N+ 3 |x| p_5/\pi.$$
This term is always non-positive by the inequality $|x|\le a(N)$ and the
explicit formula~(\ref{eqn:circ-perim}).  It is zero iff $x=0$. Hence, a
counterexample remains a counterexample even after the edge deletion
process. We drop the primes in notation, and consider
Inequality~\ref{eqn:CE} as applying to a partition of the sphere that
has no areas less than $a(N)$.

It is shown in \cite{H} that configuration can be modified so that it
is still a counterexample, but so that every region is connected and simply
connected, and so that the degree at every marked point is $3$.  We
make this modification here too.

Inequality~\ref{eqn:CE} can be rewritten as
    $$\sum_P \left[(L(P)  - T(P) B' + (5-N(P)) p_5' - \rho(P) p_5 \right] <0.$$
This uses the Euler relation
    $$\sum_P (N(P)-5) = 0$$
and the skew symmetry of $\tau$:
    $$\sum_P T(P) = 0.$$
In this form, it is clear that the inequality contradicts the
isoperimetric inequality for pentagons. This proves that there is no
perimeter strictly less than $6\,\op{reg-perim}(\pi/3,5)$.

Assume that we have a perimeter equal to $6\,\op{reg-perim}(\pi/3,5)$.
As before, we may assume that the boundaries of the regions are analytic arcs.
No edges are deleted in the edge-deletion process, because such would
lead to a strict inequality.  Thus, all areas area at least $a(N)$ and we
can apply the pentagonal isoperimetric inequalities to conclude that each
region is a regular pentagon of area $\pi/3$.  This concludes the proof.

\section{Proof of the isoperimetric inequality}

This section gives a proof of the isoperimetric inequality stated in
Section~\ref{sec:pent}.  The proof is similar to the proof of the
hexagonal isoperimetric inequality in \cite{H}.

\begin{estimate} If a region has unsigned area $x\le2\pi$, the perimeter
is at least $\op{circ-perim}(0,x)$.
\end{estimate}

\begin{proof} This is the isoperimetric inequality on a sphere, which
can be deduced from the methods of \cite{MAMM}.
\end{proof}

\begin{estimate} If an $n$-gon has unsigned area $x$ at most $2\pi$, the perimeter
is at least $\op{reg-perim}(x,n)$.
\end{estimate}

\begin{proof} This is \cite[Lemma~6.1]{HSPP}.
\end{proof}

We also use Dido's theorem for a unit sphere.

\begin{thm} \label{thm:dido} (Dido)  Let $R$ be
a region of area $x$ at most $\pi$ on the unit sphere
bounded by a great circle and some other
curve $C$.  The length of $C$ is at least that of a semicircle meeting
the great circle at right angles and enclosing an area $x$.
\end{thm}

That is, the length is at least $\op{circ-perim}(0,2x)/2$.

\begin{proof}  This follows easily from the isoperimetric inequality.
\end{proof}

\begin{lemma} If $0\le x \le y< 4\pi$,  then
    $$y\,\op{circ-perim}(0,x) \ge x\, \op{circ-perim}(0,y).$$
\end{lemma}

\begin{proof} This is an elementary consequence of the explicit
Formula~(\ref{eqn:circ-perim}).
\end{proof}

\begin{lemma} If $x\in[0.03,2.33]$ and $n\in[3,4,5,6,7]$, then
    $$\partial_1^2 \op{reg-perim}(x,n)<0.$$
\end{lemma}

\begin{proof}
This is an interval arithmetic calculation based on the explicit formula
for $\op{reg-perim}$.  If $f$ is a function, the Mathematica procedure
to check that its second derivative is negative on $[a,a+n w]$ is three
lines of code (listed below).  The result is thus readily checked (with
$w=1/1000$).
\end{proof}

\bigskip
\begin{alltt}
maxSecond[f_,a_,n_,w_]:= Module[\{i,der2\},
    der2 = D[f[x],\{x,2\}];
    Table[der2/.\{x->Interval[\{a+i w,a+(i+1)w\}]\},\{i,0,n-1\}]//Max
    ];
\end{alltt}
\bigskip

Let $$I(\ell,n,t,\rho) =\ell - t B' + (5-n) p_5' - \rho p_5.$$ To prove
the pentagonal isoperimetric inequality, we wish to show that
$I(L(P),N(P),T(P),\rho(P))$ is non-negative. Let $x_i = x(e_i)$ and $t_i
= \tau(x_i)$, and $T_{\op{abs}} = \sum |t_i|$. Let $\ell_i = L(e_i)$.
Recall that $|t_i|\le \tau_0$.

We consider several cases.  In each successive case, we may assume that
at least one of the defining conditions of the previous cases fail. A
common technique will be to reflect one arc of the perimeter $\Gamma$ in
such a way that the area increases (but not beyond $2\pi$) and the
perimeter remains the same length.  Using the negativity of the second
derivative of $\op{circ-perim}$ and $\op{reg-perim}$, most of the cases
reduce to checking that the inequality is satisfied at the endpoints of
its domain.

\subsection{Case $n\ge5$ and $T_{\op{abs}} \ge 1.301$.}  Recall that we have
chosen $x_i$ so that $|x_i|\le\pi$, so that Dido's theorem applies.
In this case, the
inequality follows because the lemma gives
    $$\ell_i\ge|t_i|\op{circ-perim}(0,2\tau_0)/(2\tau_0)$$ so that
    $$I \ge I(5,1.301\,\op{circ-perim}(0,2\tau_0)/(2\tau_0),-1.301,1) >0.$$
We now assume that $T_{\op{abs}} < 1.301$.

\subsection{Case $n\ge5$ and $\op{area(P)}\ge3.474$}
By the previous case, $|T(P)|\le T_{\op{abs}} < 1.301$.  Then
    $$I\ge I(5,\op{circ-perim}(0,3.474),-1.301,1) > 0.$$
From now on, we assume that if $n\ge5$, then we have $\op{area(P)}<3.474$.

\subsection{Case $n\ge5$ and some $x_i >\tau_0$}

If some $x_i > \tau_0$ and some $x_j < -\tau_0$, we may shrink both
bounding curves to decrease $|x_i|$, $|x_j|$, and $L(P)$, while fixing
$T(P)$ and $\op{area}(P)$.  This decreases $I$ and transforms any
counterexample to a counterexample.  If some $x_j< -\tau_0$ we may
shrink the bounding curve $L(P)$ to increase $x_j$ and $\op{area}(P)$,
while maintaining $T(P)$.  This again transforms counterexamples to
counterexamples.  If $\op{area}(P)$ ever increases to $3.474$, we are
done by the earlier case.  Thus, we may assume that all $x_j
\ge-\tau_0$.

Assume some $x_i > \tau_0$.   Let $t_{neg} = \sum_i \min(0,t_i)$.
We have $T(P)\ge \tau_0 + t_{neg}$.  Since $x_i\ge-\tau_i$, we have
$\sum_i \min(0,x_i) = t_{neg}$.  If we reflect each negative area $x_i$
about the segment $f_i$ then $P$ is transformed into a region $P'$ of
area $\op{area}(P)+2|t_{neg}|\ge \rho(P)\pi/3 + 2|t_{neg}|$.  The area of $P'$
is at most $3.474+ 2(1.301)< 2\pi$, so that it falls within the monotonic
range of $\op{circ-perim}$.  We have
    $$I \ge I(5,\op{circ-perim}(0,\rho\pi/3+2|t_{neg}|),\tau_0-|t_{neg}|,\rho).$$
The function on the right is positive for all $\rho\in[0,1]$ and all
$|t_{neg}|\in[0,1.301]$, as the explicit formula~\ref{eqn:circ-perim} readily shows.

From now on, we assume that if $n\ge5$, then we have $x_i\le\tau_0$ for
all $i$.  Thus, $x_i = t_i$.  If some $t_i<0$ and another $t_j>0$, we may
shrink the bounding curve while maintaining $t_i+t_j$ and area.  This
transformation takes counterexamples to counterexamples.  Thus, we may
assume that all $t_i$ have the same sign.

\subsection{Case $n= 5$, $\rho\ge0.9957$, $T(P)\in[-0.0711,0.117]$,
    and some edge $f_i$ has arc-length at least $1$.}

\begin{figure}[htb]
  \centering
  \includegraphics{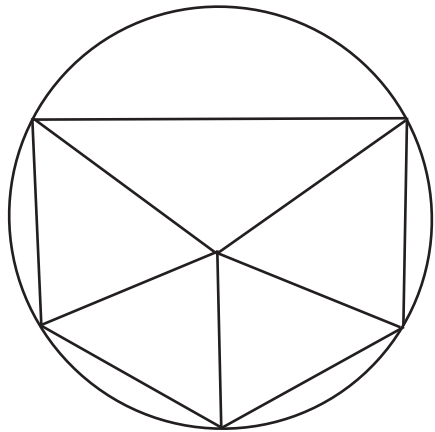}
  \caption{Perimeter Minimizer with a long edge.  Four edges have length $u$
    and the top edge has length $1$.  The angles at the origin are
    $\alpha$, $\beta$, $\beta$, $\beta$, $\beta$.  The arclength of the radius
    is $x$.  Also, the angle $\gamma$ is that along the isosceles triangle
    with base $1$ and $\delta$ is that along the isosceles triangle with
    base $u$.}
  \label{diag:inscribed}
\end{figure}

\begin{lemma} Let $P$ be a spherical pentagon of area $x < 2\pi$ on the unit
sphere.  Assume that it has an edge of arclength at least $1$.  If the
regular pentagon of area x has edges of arclength at least $1$, then that
pentagon is perimeter minimizing.  Otherwise the perimeter minimizer is
that with edges $1$, $u,u,u,u$ inscribed in a circle.
(See Figure~\ref{diag:inscribed}.)
\end{lemma}

\begin{proof} The case of $5$ equal edge lengths has already been
considered.  Let us assume that the edge constraint rules out the
regular pentagon.  The proof of \cite[Lemma~6.1]{HSPP} shows that the
optimal solution will have edge $1$, $u$, $u$, $u$, $u$, for some $u$.

The argument that the optimal figure will have a circumscribing circle
is classical, at least in the analogous case of the plane. In \cite{HW},
it is shown that the optimal polygon among equilateral polygons is the
regular polygon.  Their argument actually shows that for any three
consecutive equilateral sides, the two interior angles along the middle
segment are equal.  This implies that the pentagon will have equal
angles at vertices at the end of two edges of length $u$.  For a
pentagon, this is sufficient to give a circumscribing circle.
\end{proof}

By simple spherical trig (the law of cosines), the perimeter and area
of this pentagon are described parametrically in terms of the arclength $x$
of the radius of the circle.
This gives the following parametric equations for perimeter and area.
%
\bigskip
\begin{alltt}
perimArea[x_]:= Module[\{beta,alpha,u,gamma,delta,perim,area\},
    alpha = ArcCos[(Cos[1] - Cos[x]^2)/(Sin[x]^2)];
    beta = (2 Pi - alpha)/4;
    u = ArcCos[Sin[x]^2 Cos[beta] + Cos[x]^2];
    gamma = ArcCos[(Cos[x]-Cos[x] Cos[1])/(Sin[x] Sin[1])];
    delta = ArcCos[(Cos[x]-Cos[x] Cos[u])/(Sin[x] Sin[u])];
    perim = 1 + 4 u;
    area = (alpha + 2 gamma - Pi) + 4 (beta+ 2delta - Pi);
    \{perim,area\}
    ];
\end{alltt}
\bigskip

The signed area bounded by the generalized polygon with segments $f_i$
is at least $\rho\pi/3 - T(P)$, so that $T(P)\ge \rho\pi/3-\op{area}_x$.
This gives
    \begin{equation}
    I \ge I (5,\op{perim}_x,\rho\pi/3-\op{area}_x,\rho)
    \ge I(5,\op{perim}_x,\pi/3-\op{area}_x,1)
    \label{eqn:parametric}
    \end{equation}
The perimeter is increasing in $x$, for $x\le\pi/2$.  Clearly $x\ge0.5$.
If $x>0.7$, then $\op{perim}_x>3.8$ and
    $$I\ge I(5,3.8,-0.0711,1)>0.$$
The area is increasing in $x$, for $x\le0.7$.  If $x<0.6$, then
    $$0.92 < 0.9957\pi/3 - 0.117\le \rho\pi/3 - T(P)\le\op{area}_x < 0.82.$$
This contradiction shows that we may assume $x\in[0.6,0.7]$.  For
$x\in[0.6,0.7]$, the right-hand side of Inequality~\ref{eqn:parametric}
is positive, as an easy calculation shows.

\subsection{Case $n=5$, $\rho\ge0.9957$, $T(P)\in[-0.0711,0.117]$}
\label{sec:C}

This subsection follows the proof of the hexagonal honeycomb conjecture
closely.  See \cite[Sec. 8]{H}. By the isoperimetric inequality, we may
assume that all edges $e_i$ are arcs of circles.  We may assume that
each arc has the same curvature by \cite[Thm 2.3,~step~3]{M94}.
Following \cite[Prop.~6.1]{H}, the perimeter is at least
    $$
    \begin{array}{lll}
    &\op{circ-perim}(1,|t|/\op{reg-perim}(\rho\pi/3-t,5))
        \op{reg-perim}(\rho\pi/3-t,5)\\
        &\qquad\qquad \le \op{circ-perim}(1,|t|/3.65)\op{reg-perim}(\rho\pi/3-t,5).
    \end{array}
    $$
Let $C(\rho,t)$ be the right-hand side of this inequality. We have
$$I \ge I(5,C(\rho,t),t,\rho).$$

The right-hand side is non-negative.  We first use interval arithmetic
to prove positivity outside the set
    $$\{(\rho,t) : \rho\in[0.999,1]\quad t\in[-0.04,0.043]\}.$$
Then the $\rho$-derivative is shown to be negative by another interval
arithmetic calculation, and this reduces the problem to $\rho=1$. Yet
another interval calculation reduces the problem to the interval
    $$t\in[-0.005,0.005].$$  Finally an interval calculation shows that
the second derivative with respect to $t$ is strictly positive for
$t\in[-0.005,0.005]$.  All of these calculations were made with
Mathematica's interval function and a few lines of code.  Exact
arithmetic shows that the function is zero and has zero derivative at
$t=0$. Thus, $t=0$ is the unique minimizer. Equality is achieved iff
$\rho=1$ and $t=0$. In this case, we have a regular pentagon of area
$\pi/3$.

From here on, we assume that if $n=5$, we have $\rho\le0.9957$ or
    $T(P)\not\in[-0.0711,0.117]$.

\subsection{Case $n\ge5$, $T(P)\ge0$}  Define $\rho_1$ by
    $$
    \rho_1(n)=\begin{cases}
        0.913 & n=5\\
        0.952 & n=6\\
        0.99& n=7\\
        1.0 & n\ge8.
        \end{cases}
    $$
Define $T_1$ by
    $$
    T_1(n) = \begin{cases}
        0.117 & n=5\\
        0.065 & n=6\\
        0.0134 & n=7
        \end{cases}
    $$

For $\rho\le\rho_1(n)$, we find that
$$I \ge I(n,\op{circ-perim}(0,\rho\pi/3),0,\rho) >0$$
by the explicit formula for $\op{circ-perim}$.

Assume $\rho>\rho_1(n)$.  For $T(P)\ge T_1(N(P))$ we have
    $$I \ge I(n,\op{circ-perim}(0,\rho\pi/3),t,\rho) >0,$$
as we can readily check by looking at the endpoints of functions with
negative second derivatives.

Assume $\rho\ge\rho_1(n)$ and $T(P)\in[0,T_1(N(P))]$.  We have reduced
to the case $N(P)\in[5,6,7]$.  If $N(P)=5$, we assume that
$\rho\le0.9957$. We have
$$I \ge I(n,\op{reg-perim}(\rho\pi/3-t,n),t,\rho)>0$$
by checking endpoints and using a second derivative test.

\subsection{Case $n\ge5$, $T(P)\le0$}  Define $T_0$ by
    $$
    T_0(n) = \begin{cases}
        -0.1382 & n=5\\
        -0.0711 & n=6\\
        -0.01362& n=7\\
        0 & n\ge8.
        \end{cases}
    $$
Assume that $T(P)\in[-1.301,T_0(n)]$.  We have
    $$I\ge I(n,\op{circ-perim}(0,\rho\pi/3-2t),t,\rho)>0$$
by a second-derivative test.

Assume that $T(P)\in[T_0(n),0]$. We may also assume $n<8$.   Define
$\rho_2$ by
    $$
    \rho_2(n) = \begin{cases}
        0.913 & n=5\\
        0.952 & n=6\\
        0.99  & n=7.
        \end{cases}
        $$
If we also assume $\rho\in[0,\rho_2(N(P))]$, we find that
    $$
    I\ge I(n,\op{circ-perim}(0,\rho\pi/3-2t),t,\rho)>0.
    $$

Assuming that $\rho\in[\rho_2(n), 0.9957]$ when $n=5$, and
$\rho\in[\rho_2(n),1]$ otherwise, we have
$$I \ge I(n,\op{reg-perim}(\rho\pi/3-t,n),t,\rho) >0$$
for all $T\in [T_0(n),0]$ by a second derivative test.


We have now completed the proof of the pentagonal isoperimetric inequality
for $n\ge5$.

\subsection{Case $n=4$}
\label{sec:n=4}

This case is similar to the case $n\ge5$.  We simply list the
differences. Note that $T(P)=\sum^4 t_i\ge -4\tau_0$. If $T(P)>0.168$,
use
$$I \ge I(4,\op{circ-perim}(0,\rho\pi/3),0.168,\rho)>0.$$
If some area $|x_i|$, $|A_P|$, or $|\op{area}(P)|$ has minimal
real representative between $3.56$ and $4\pi-3.56$, then
    \begin{equation}
    I \ge I(4,\op{circ-perim}(0,3.56),-4\tau_0,1)>0.
    \label{eqn:4tau0}
    \end{equation}

If some $x_i>\tau_0$, we contract the boundary as before to assume that
$x_j\ge-\tau_0$ for all $i$.  Reflecting the region corresponding to
each $t_j<0$, we add $2|t_j|$ to the area. The area increases by at most
$2\tau_0$ to at most $3.56+2\tau_0$.  If this is greater than $3.56$, it
is less than $2\pi$, and the area estimate~(\ref{eqn:4tau0}) applies.
Thus, we may assume after each reflection that the new area is less than
$3.56$. We obtain
$$I \ge I(4,\op{circ-perim}(0,\rho\pi/3+2|t_{neg}|,\tau_0-|t_{neg}|,\rho).$$
The right-hand-side is positive for $\rho\in[0,1]$ and $|t_{neg}|\in[0,3\tau_0]$.

As above, we may assume that $x_i\in[-\tau_0,\tau_0]$.  Recall that
$a(4)=0.1$.  We  finish the case $n=4$ with several estimates.
    $$I\ge I(4,\op{circ-perim}(\rho\pi/3),t,\rho) >0,$$
if $\rho\in[0.1,0.85]$ and $t\in[0,0.168]$.
    $$I\ge I(4,\op{reg-perim}(\rho\pi/3-t,4),t,\rho)>0,$$
if $\rho\in[0.85,1]$ and $t\in[0,0.168]$.
    $$I\ge I(4,\op{circ-perim}(\rho\pi/3-2t),t,\rho)>0$$
if $\rho\in[0.1,0.85]$ and $t\in[-4\tau_0,0]$, or if
    $\rho\in[0.1,1]$ and $t\in[-4\tau_0,-0.25]$.
    $$
    I\ge I(4,\op{reg-perim}(\rho\pi/3-t,4),t,\rho)>0$$
if $\rho\in[0.85,1]$ and $t\in[-0.25,0]$.

\subsection{Case $n=3$}

This case is similar to $n=4$. If $T(P)>0.22$, we use an estimate like
that of Section \ref{sec:n=4}. The proof is the same if we modify the
constants in the proof as follows:
    $$
    \begin{array}{lll}
    -0.25\mapsto-0.35,& 3.56\mapsto 2.84,&
        0.168\mapsto 0.22,\\
     4\mapsto 3,& 3\mapsto 2,
        & 0.85\mapsto 0.8.
    \end{array}
    $$

\subsection{Case $n=2$, Great Circle}

\begin{figure}[htb]
  \centering
  \includegraphics{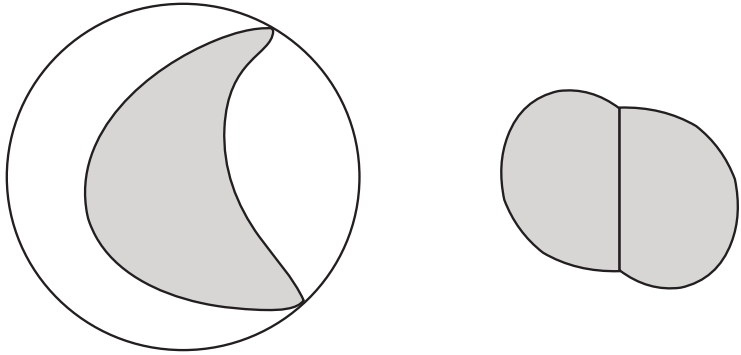}
  \caption{Two types of digons.}
  \label{diag:digons}
\end{figure}

There are two types of digons.  We call them {\it great circle digons}
or {\it simple arc digons}.  Recall that we choose $x_i$ so that its
minimal real representative lies in the interval
    $\leftclosed -\pi,\pi \rightopen$.  A choice of arc $f_i$ of a great
    circle was chosen for each edge $e_i$ to make the minimal real
    representatives lie in the given interval.  In the case of a digon
    the two arcs $f_1$ and $f_2$ have the same endpoints.  If the endpoints
    are the two parts of a great circle then we have a great circle
    digon.  If the arcs trace the same path with opposite orientations then
    we have a simple arc digon.
These two types are shown in Figure~\ref{diag:digons}.

Assume we have a great circle digon.  We have $|A_P| = 2\pi$,
$x_1,x_2\in[-\pi,0]$, and $|x_1|+|x_2|+|\op{area}(P)|=2\pi$.

If the area of the digon is at
least $2.29$, we have
    \begin{equation}
    I \ge I(2,\op{circ-perim}(0,2.29),-2\tau_0,1)>0.
    \label{eqn:digon}
    \end{equation}

\begin{figure}[htb]
  \centering
  \includegraphics{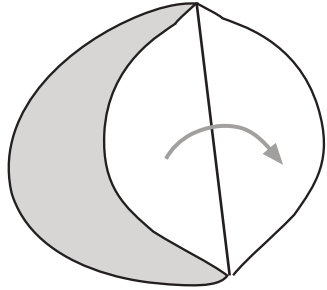}
  \caption{Reflection on a digon region.}
  \label{diag:reflect}
\end{figure}

We claim that if the unsigned area $|\op{area}(P)|$ of the digon  is in
$[\pi/3,2.29]$, then the area increases to a value in
$[2.29,4\pi-2.29]$, upon reflection of the region represented by $x_1$
(Figure~\ref{diag:reflect}). This implies that $I>0$ by
Inequality~\ref{eqn:digon}.
 In fact, after reflection,
the area is at most
    $$|\op{area(P)}| + 2|x_1|\le 2.29+2\pi<4\pi-2.29.$$
Recall that
    $$|x_1| = 2\pi-|\op{area}(P)|-|x_2|\ge 2\pi-2.29-\pi=\pi-2.29.$$
Hence, the area is at least
    $$|\op{area}(P)|+2|x_1|\ge \pi/3 + 2 (\pi-2.29) > 2.29.$$

Finally, suppose that the unsigned area of the digon is at most $\pi/3$.
In this case, we again reflect the areas represented by $x_1$.
We have
    $$|x_1| = 2\pi-|\op{area}(P)| -|x_2|\ge 2\pi-\pi/3-\pi = 2\pi/3.$$
So the reflected area $R_A$ satisfies
    $$2.29 <2|x_1| \le R_A \le \pi/3 + 2\pi < 4\pi - 2.29.$$
Hence, we may again apply Inequality~\ref{eqn:digon}.  This completes
the case of great circle digons.

\subsection{Case $n=2$, Simple Arc}

In this case, $A_P=0$ and the signed area is simply $x_1+x_2$.

If $T(P)\ge0.271$, then the result follows from
$$I \ge I(2,\op{circ-perim}(0,\rho\pi/3),0.271,\rho)>0.$$
Assume that $T(P)\le 0.271$.

If the area is between $2.29$ and $2\pi$, we may apply Inequality~\ref{eqn:digon}
to prove that $I>0$.  Assume that the area is less than $2.29$.

If some $x_i < -\tau_0$ contract the perimeter increasing $x_i$ and
the area $x_1+x_2$ until either the area becomes $2.29$ or $x_i\ge-\tau_0$.
Assume now that $x_i\ge-\tau_0$.

If some $x_i > \tau_0$ contract the perimeter while decreasing $x_1+x_2$
until either $x_i=\tau_0$ or $x_1+x_2\le\pi/3$.  Assume $x_1+x_2\le\pi/3$.
We have
    $$T(P) = \tau_0 + t_{\min}\le 0.271$$
and
    $$x_{\min}\le t_{\min} \le 0.271-\tau_0 = -0.049.$$
Reflecting $x_{\min}$, we have
    $$I\ge I(2,\op{circ-perim}(0,\rho\pi/3+2|t_{\min}|),\tau_0-|t_{\min}|,\rho)
        >0,$$
for $|t_{\min}|\in [0.049,\tau_0]$ and $\rho\in[0,1]$.

If $x_1,x_2\in[-\tau_0,\tau_0]$ and they have opposite signs,
we may contract the perimeter while
maintaining $x_1+x_2=t_1+t_2$, until they have the same sign.

Finally, assume that $x_1,x_2\in[-\tau_0,\tau_0]$ and that they have
the same sign.
    $$0.271 \ge T(P) = t_1+t_2 = x_1+x_2 = |\op{area}(P)|\ge0,$$
so the sign is positive.  Assume that $\rho\ge a(2)=0.1$ (as allowed by
the hypothesis of the pentagonal isoperimetric inequality).  We have
    $$I \ge I(2,\op{circ-perim}(0,\rho\pi/3),\rho\pi/3,\rho) >0,$$
for $\rho\in[0.1,1]$.

\section{Concluding Remarks}

The method of the honeycomb conjecture is seen to extend to the
spherical case.  As a further test of the generality of the method, it
would be interesting to see whether the method can be used to prove the
corresponding results for $N=3$, $4$, and $6$ regions of equal area. A
further test of the method would be to extend it to the honeycomb
problem in the hyperbolic plane.  I am not aware of any obstructions to
doing so.

My motivation for pursuing the case of a sphere is that for
    $$N\ne1,2,3,4,6,12,$$
the optimal solution will not consist of congruent regular polygons.  The
study of these irregular cases can lead us to valuable insights into
problems such as the Kelvin problem, where the best known solution contains
two types of cells.  As a preliminary step, it is necessary to develop
the spherical theory at a level comparable to what has been done for the
plane.  That is the purpose of these calculations.



\begin{thebibliography}{99}

\bibitem{A} F. J. Almgren, Jr.  Existence and regularity of almost everywhere
    of solutions to elliptic variational problems with constraints,
    Mem. AMS, 165 (1976).

\bibitem{C} C. Cox, L. Harrison, M. Hutchings, S. Kim, J. Light, A. Mauer, M.
    Tilton, The shortest enclosure of three connected areas in $\ring{R}^2$,
    Real Analysis Exchange, Vol. 20(1), 1994/95, 313--335.

\bibitem{F} H. Federer, Geometric Measure Theory, Springer-Verlag, 1969.

\bibitem{FT43} L. Fejes T\'oth, \"Uber das k\"urzeste Kurvennetz
    das eine Kugeloberfl\"ache in fl\"achengleiche konvexe
    Teil zerlegt, Mat. Term.-tud. \'Ertesit\"o 62 (1943), 349--354.

\bibitem{FT64a} L. Fejes T\'oth, Regular Figures, MacMillan Company, 1964.

\bibitem{FT64b} L. Fejes T\'oth, What the bees know and what they do not know,
    Bulletin AMS, Vol 70, 1964.

\bibitem{HW} W. Habicht and B. L. van der Waerden, Lagerung von
    Punkten auf der Kugel, {\it Math. Ann.} {\bf 123} (1951) 223-234.

\bibitem{HSPP} T. C. Hales, The sphere packing problem, J. Comp. and
Appl. Math. 44 (1992) 41--76.

\bibitem{H} T. C. Hales, The  honeycomb conjecture,
    Discr. Comput. Geom. 25:1-22 (2001).


\bibitem{M94} F. Morgan, Soap bubbles in ${\ring{R}^2}$ and in surfaces,
    Pacific J. Math, 165 (1994), no. 2, 347--361.

\bibitem{M95} F. Morgan, Geometric Measure Theory, A Beginner's Guide, Second
Edition, Academic Press, 1995.

\bibitem{M99} F. Morgan, The hexagonal honeycomb conjecture,
    Trans. AMS, Vol 351, Number 5, pages 1753--1763,
    1999.

\bibitem{MAMM} F. Morgan, The Isoperimetric Problem on Surfaces,
    Amer. Math. Monthly, vol 106, 5, May 1999, 430--439.

\bibitem{T} J. Taylor, The structure of singularities in soap-bubble-like
    and soap-film-like minimal surfaces, Annals of Math., 103 (1976),
    489-539.


\end{thebibliography}
\end{document}